# Minimal Leader Set for Controllability of $k$-distant Trees

Li Dai[a]*

[a] *Department of Mathematics, National University of Defense Technology, Changsha, China*
*e-mail: daili@nudt.edu.cn

Abstract:

Minimal controllability problem plays an important role in the field of network control. A New concept—Minimum Perfect Critical Set (MPCS) is proposed. Four different MPCSs were found for $k$-distant tree graphs. Based on this concept of MPCS, an algorithm for finding the minimal leader set is provided. Numerical experiments show that these theories enable the algorithm to find a minimal leader set with a probability of more than 0.98. Further, some other numerical characteristics of the minimal leader set of $k$-distant trees were found.

## 1. Introduction

Due to its wide application, people pay more and more attention to the controllability of undirected networks. For example, Figure 1 is a typical UAV swarm architecture [1]. The air-air communication links can be modeled as an undirected network. In this kind of system, individuals controlled by the ground are called leaders. In contrast to the various insights and results on the controllability of directed networks, there are still many basic unsolved problems in the study of controllability of undirected networks. To ensure the controllability of an undirected network, it is necessary to determine how to select the leader(leaders), what is the minimum number of leaders, how to get the minimum set of leaders, how many minimum sets of leaders are there in the network. These problems are interesting because of their potential applications to various fields such as complex networks, multi-agent systems, biology, etc.

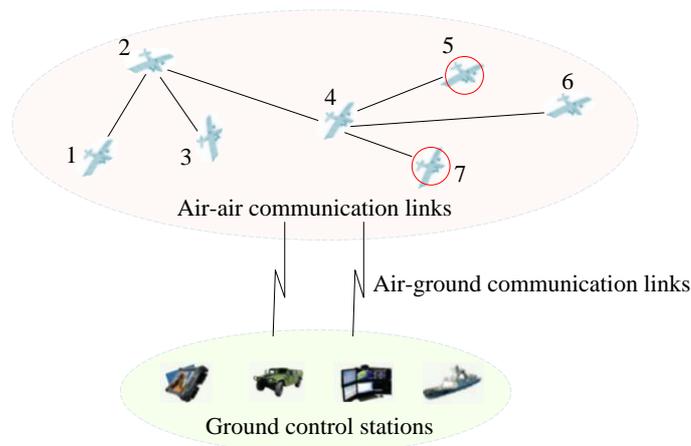

Figure 1  The UAV swarm architecture

## 1.1 Graphical Model and Notations

Regard the UAVs as the vertices of the graph, the air-air communication links between the UAVs as the corresponding edges, we will get the graphical model of a UAV communication network.

Let $G = (V, E)$ be an undirected and unweighted simple graph, where $V = \{v_1, v_2, \cdots, v_n\}$ is a vertex set and $E = \{v_i v_j | v_i \text{ and } v_j \in V\}$ is an edge set, where an edge $v_i v_j$ is an unordered pair of distinct vertices in $V$. If $v_i v_j \in E$, then $v_i$ and $v_j$ are said to b *adjacent* or *neighbors*. $N_S(v_i) = \{v_j \in S | v_i v_j \in E(G)\}$ represents the neighboring set in $S$ of $v_i$, where $S \subset V$. The cardinality of $S$ is denoted by $|S|$. $d_G(v)$ is the degree of vertex $v$ in $G$ and $d_G(v) = |N_G(v)|$. If $d_G(v) = 1$, then $v$ is called a pendent vertex. $G[S]$ is the induced subgraph, whose vertex set is $S$ and edge set is $\{v_i v_j \in E(G) | v_i, v_j \in S\}$. The *valency matrix* $\Delta(G)$ of graph $G$ is a diagonal matrix with rows and columns indexed by $V$, in which the $(i,i)$-entry is the degree of vertex $v_i$, e.g., $|N_G(v_i)|$. Any undirected simple graph can be represented by its *adjacency matrix*, $D(G)$, which is a symmetric matrix with 0-1 elements. The element in position $(i,j)$ in $D(G)$ is 1 if vertices $v_i$ and $v_j$ are adjacent and 0 otherwise. The symmetric matrix defined as $\boldsymbol{L}(G) = \boldsymbol{\Delta}(G) - \boldsymbol{D}(G)$ is the *Laplacian* of $G$. Figure 2 is the undirected graphical model of the UAV communication network in Figure 1 and the Laplacian of $G$ is

$$\boldsymbol{L}(G) = \begin{bmatrix} 1 & -1 & 0 & 0 & 0 & 0 & 0 \\ -1 & 3 & -1 & -1 & 0 & 0 & 0 \\ 0 & -1 & 1 & 0 & 0 & 0 & 0 \\ 0 & -1 & 0 & 4 & -1 & -1 & -1 \\ 0 & 0 & 0 & -1 & 1 & 0 & 0 \\ 0 & 0 & 0 & -1 & 0 & 1 & 0 \\ 0 & 0 & 0 & -1 & 0 & 0 & 1 \end{bmatrix}.$$

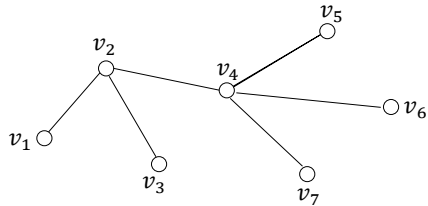

Figure 2  undirected graphical model $G$

## 1.2 Literature Review

The neighbor-based controllability of the undirected graph under a single leader was formulated by Tanner in [2] and a necessary and sufficient condition expressed in terms of eigenvalue and eigenvector was derived. In the case of multiple leaders, the problem of finding a minimal set of leaders was considered and proved to be NP-hard in [3] and [4]. Some algebraic conditions were developed [5,6,7,8] and these algebraic conditions lay the foundation for understanding interaction between topological structures of an undirected graph and its controllability. They also served as the theoretical basis of this paper. The research efforts on characterizing the controllability from a graphical point of view were also motivated by [2] to build controllable topologies. Many kinds of noncontrollable topologies were characterized, such as asymmetric graph with respect to the anchored nodes [9], quotient graphs [10], nodes with the same number of neighbors [11], controllability destructive nodes [8], etc. Useful tools and methods were developed to study the controllability of undirected graph, such as downer branch for tree graphs [12], Zero forcing set, balancing set and color change rule for weighted undirected graph [13], equitable partitions [5,14,15,16,17,18], leader and follower subgraphs [11], $\lambda$-core vertex [19,20], NEPSes of graphs[21], etc. Omnicontrollable systems are defined by [20], in such systems, the choice of leader vertices that control the follower graph is arbitrary. In the study [22], two algorithms are established for selecting the fewest leaders to preserve the controllability and the algorithm for leaders' locations to maximize non-fragility is also designed.

Although many scholars have achieved many remarkably strong and elegant results, this problem has not been solved yet. As it is well known that any undirected simple connected graph on $n$ vertices is always $(n-1)$-omnicontrollable, the choice of leader vertices for minimal controllability is important.

Therefore, our aim is to find a method for giving a direct interpretation of the leader vertices from a graph-theoretic vantage point. In this sense, we provide a new concept, perfect critical vertex set, to identify the potential leader vertices.

## 1.3 Preliminary Results

Throughout this paper, it is assumed without loss of generality that $F$ denotes a follower set; its vertices play the role of followers and the vertices in $\overline{F}$ are leaders(driver nodes), where $\overline{F} = V \setminus F$ denotes the complement set of $F$. Let $\boldsymbol{y}$ be a vector and $\boldsymbol{y}_S$ denote the vector obtained from $\boldsymbol{y}$ after deleting the elements in $\overline{S}$. Let $\boldsymbol{L}_{S \to T}$ denote the matrix obtained from $\boldsymbol{L}$ after deleting the rows in $\overline{S}$ and columns in $\overline{T}$. If the followers' dynamics is (see [2])

$$\boldsymbol{x}' = \boldsymbol{A}\boldsymbol{x} + \boldsymbol{B}\boldsymbol{u},$$

where $\boldsymbol{x}$ is the vector of all $x_i$ corresponding to follower $v_i \in F$ and $\boldsymbol{u}$ is the external control input vector that is injected to only leaders, then, the system is controllable if and only if the $N \times NM$ controllability matrix

$$\boldsymbol{C} = [\boldsymbol{B}, \boldsymbol{AB}, \boldsymbol{A}^2\boldsymbol{B}, \cdots, \boldsymbol{A}^{N-1}\boldsymbol{B}]$$

has full row rank, where $\boldsymbol{B} = \boldsymbol{L}_{F \to \overline{F}}$ and $\boldsymbol{A} = \boldsymbol{L}_{F \to F}$. This mathematical condition for controllability is well known as Kalman's controllability rank condition [23,24,25].

For example, see Figure 2, let $\overline{F} = \{v_1, v_4, v_6\}$, then $\boldsymbol{A} = \begin{bmatrix} 3 & -1 & 0 & 0 \\ -1 & 1 & 0 & 0 \\ 0 & 0 & 1 & 0 \\ 0 & 0 & 0 & 1 \end{bmatrix}$, $\boldsymbol{B} = \begin{bmatrix} -1 & -1 & 0 \\ 0 & 0 & 0 \\ 0 & -1 & 0 \\ 0 & -1 & 0 \end{bmatrix}$, $rank(\boldsymbol{C}) = 3 < 4$, the network is NOT controllable under the leader set $\overline{F}$. However, if $v_5$ is selected as the leader vertex, e.g., $\overline{F} = \{v_1, v_5, v_6\}$, the system is controllable. Therefore, it is very important to select the leaders properly for the network to be controllable.

In terms of eigenvalues and eigenvectors of Laplacian submatrices, [5,8] presented a necessary and sufficient algebraic condition on controllability.

**Proposition 1**[5,8]

The undirected graph $G$ is controllable under the leader set $\overline{F}$ if and only if $\boldsymbol{y}_{\overline{F}} \neq \boldsymbol{0}$ ($\forall \boldsymbol{y}$ is an eigenvector of $\boldsymbol{L}$).

The remainder of this paper is organized as follows. Section 2 presents the new concept of MPCS. In addition, it presents the necessary and sufficient conditions for $S$ to be a minimum perfect critical 2 set and proves an interesting result that a minimum perfect critical 3 set does not exist. Section 3 constitutes the main part of this paper. It presents some special MPCS of $k$-distant graphs. Section 4 presents three examples to illustrate how to use the theories proposed to find MPCSs. Finally, Section 5 concludes the paper.

## 2. Minimal Perfect Critical Set

### 2.1 Three New Definitions

According to Proposition 1, for any $S \subset V$ and $S \neq \emptyset$, if there exists an eigenvector $\mathbf{y}$ of Laplacian matrix $\mathbf{L}$ such that $\mathbf{y}_S = \mathbf{0}$, then $S$ cannot be used as a leader set. Therefore, to locate the leaders of graph $G$, the following concepts are proposed.

**Definition 1**

Let $S$ be a nonempty subset of $V$. If there exists an eigenvector $\mathbf{y}$ such that $\mathbf{y}_{\bar{S}} = \mathbf{0}$, then $S$ is called a *critical set*(CS) and $\mathbf{y}$ is an inducing eigenvector. $S$ is called a $k$ critical set if $|S| = k$.

**Definition 2**

Let $S$ be a critical set. If there exists an eigenvector $\mathbf{y}$ satisfying $\mathbf{y}_{v_i} \neq \mathbf{0}$ ($\forall v_i \in S$), then $S$ is called a *perfect critical set* (PCS). $S$ is called a $k$ perfect critical set if $|S| = k$.

**Definition 3**

A perfect critical set is called a *minimum perfect critical set* (MPCS) if any proper subset of it is no longer a PCS. $S$ is called a $k$ minimum perfect critical set if $|S| = k$.

From the definitions stated above, $V$ is a trivially CS or PCS induced by the eigenvector $\mathbf{1}_n$. $V$ is an MPCS if and only if $G$ is controllable under any single vertex selected as a leader, e.g., $G$ is omnicontrollable.

For example, see Figure 2, $\{v_3, v_5\}$ is not a CS; $\{v_1, v_2, v_3, v_5, v_6, v_7\}$ is a CS but not a PCS. $\{v_1, v_3, v_5, v_6, v_7\}$ is a PCS but not a MPCS. $S_1 = \{v_1, v_3\}, S_2 = \{v_5, v_6\}, S_3 = \{v_5, v_7\}, S_4 = \{v_6, v_7\}$ is all of the MPCSs of the graph $G$ in Figure 2.

### 2.2 New way to find the minimum leader set

MPCS defined above will provide a new way to find the minimum leader set. Recall Proposition 1 and the definition of MPCS, we can restate Proposition 1 as follows:

**Remark 1**

The undirected graph $G$ is controllable under the leader set $\bar{F}$ if and only if $S \cap \bar{F} \neq \emptyset$ ($\forall S$ is a MPCS).

Let's give an example to illustrate the connections between MPCS and leaders' selection for controllability. For example, for graph $G$ in Figure 2, by Remark 1 and its MPCS stated above, the minimum leader set is $\{v_i, v_j, v_k\}$, where $v_i \in \{v_1, v_3\}$ and $v_j, v_k \in \{v_5, v_6, v_7\}$. Therefore, the minimum number of leaders is 3. In other words, when we find all MPCS of $G$, we find the minimum leader set and hence the minimum number of leader vertices. More than this, by finding out MPCSs, we will also get the number of minimum leaders set. For graph $G$ in Figure 2, there are altogether 6 different minimum leader sets.

Many MPCS have typical graphical characteristics. For example, the MPCS of $G$ in Figure 2 has the graphical structure stated in Theorem 2. Graphical structure of MPCS enables us to study the acroscopic controllability from the microscopic structure of the network. Subsection 2.3 will focus on studying the graphical structures of MPCS.

### 2.3 graphical characterization of MPCS

For an undirected graph, Laplacian matrix $\mathbf{L}$ is symmetric and all the eigenvectors are orthogonal to each other. Hence, by knowing that $\mathbf{1}_n$ is an eigenvector of $\mathbf{L}$, it is immediate that all the other eigenvectors of $L$ are orthogonal to $\mathbf{1}_n$, i.e., for all eigenvectors $y$,

$$\mathbf{1}_n^T \cdot y = \sum_{i=1}^{n} y_i = 0. \tag{1}$$

For all nontrivial graphs, if $S$ is a CS, then

$$|S| \geq 2. \tag{2}$$

In fact, suppose that $|S| = 1$; without loss of generality, $S = \{v_1\}$. Let $y = (y_1, y_2, \cdots, y_n)^T$ be the inducing eigenvector associated with eigenvalue $\lambda$. Then, $Ly = \lambda y$ and $y_{\bar{S}} = 0$. By $y_{\bar{S}} = 0$ and (1), $y_S = 0$, e.g., $y = 0$. Further, since any subset $S$ with $|S| = 1$ is not a CS, by Remark 1, $G$ is controllable with the leader set $\bar{F}$ when $|F| = 1$.

Now, we investigate the properties of a critical set. First, a sufficient condition for $S$ to be a CS is provided in the following Theorem 1, which describes a special case of the symmetry-based uncontrollability results.

Theorem 1

Let $G$ be an undirected connected graph of order $n$, $S \subset V$, and $|S| \geq 2$. If for any $v \in \bar{S}$, either $N_S(v) = \emptyset$ or $N_S(v) = S$, then $S$ is a critical set.

**Proof** Let $|\{v \in \bar{S} | N_S(v) = S\}| = m$. Then,

$$L_{S \to S} - m I_{|S|}$$

is the Laplacian of subgraph $G[S]$, where $I_{|S|}$ denotes the $|S|$-dimensional identity matrix. Considering (1), there exists an eigenvector $y_0$ of the Laplacian $L_{S \to S} - m I_{|S|}$ such that $\mathbf{1}_{|S|}^T y_0 = 0$.

Set vector $y$ as $y_S = y_0$ and $y_{\bar{S}} = 0$. It can be seen that

$$Ly = \begin{bmatrix} L_{S \to S} & L_{S \to \bar{S}} \\ L_{\bar{S} \to S} & L_{\bar{S} \to \bar{S}} \end{bmatrix} \begin{bmatrix} y_0 \\ 0 \end{bmatrix} = \begin{bmatrix} \lambda y_0 \\ L_{\bar{S} \to S} \cdot y_0 \end{bmatrix}.$$

Noting that the rows in matrix $L_{\bar{S} \to S}$ are either ones or zeros, the conclusion is proved by $\mathbf{1}_{|S|}^T y_0 = 0$. ∎

Based on the above-mentioned properties, the critical $k$ set with $k \leq 3$ can be determined directly from the graphical characterization.

This is achieved via a detailed analysis of the inducing eigenvector.

Lemma 1

Let $G$ be an undirected connected graph and $S$ be a perfect critical $k$ set. Then, for any $v \in \bar{S}$, $|N_S(v)| \neq 1$ and $|N_S(v)| \neq k - 1$.

**Proof** Let $S = \{v_1, v_2, \cdots, v_k\}$ be a perfect critical set and $y = (y_1, y_2, \cdots, y_k, 0, 0, \cdots, 0)^T$ be the inducing eigenvector. $y_i \neq 0 (\forall 1 \leq i \leq k)$ since $S$ is a perfect critical set.

$\forall v \in \bar{S}$, suppose that $|N_S(v)| = 1$. Without loss of generality, say, $vv_1 \in E(G)$ and $vv_i \notin E(G)(\forall i \neq 1)$. Then, $L_{\{v\} \to V} \cdot y = -y_1 \neq 0$. On the other hand, $y_{\bar{S}} = 0$ and $v \in \bar{S}$. Hence,

$$L_{\{v\} \to V} \cdot y = 0.$$

This is a contradiction.

Together with (1), $|N_S(v)| \neq k - 1$ can be proved similarly. ∎

By (2), a critical 2 set is also a minimum perfect critical 2 set. The following Theorem 2 will follow from Lemma 1 and Proposition 1.

Theorem 2

Let $G$ be an undirected connected graph, $S \subset V$, and $|S| = 2$. Then, $S$ is a minimum perfect critical 2 set if and only if, either $N_S(v) = \emptyset$ or $N_S(v) = S$ ($\forall v \in \overline{S}$). ∎

For example, see graph $G$ in Figure 2. All its 4 MPCSs can be recognized by the graphical characterization stated in Theorem 2.

Next, we will prove that a minimum perfect critical 3 set does not exist.

Theorem 3

Let $G$ be an undirected connected graph, $S \subset V$, and $|S| = 3$. Then, $S$ is not a minimum perfect critical set.

**Proof** Suppose that $S$ is a minimum perfect critical set. Consider the subgraph $G[S]$. All 4 possible topology structures of $G[S]$ are shown in Figure 3.

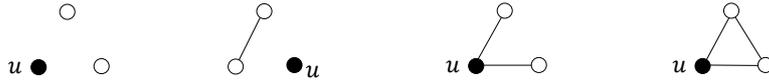

Figure 3  All 4 possible topology structures of graph with 3 vertices

For each topology of $G[S]$ in Figure 3, let $T$ be the vertex set of white vertices and $u$ be the black vertex. One can have either $N_T(u) = \emptyset$ or $N_T(u) = T$.

By Lemma 1，$\forall v \in \overline{S}$, either $|N_S(v)| = 0$ or $|N_S(v)| = 3$.

Noting that $\overline{T} = \{u\} \cup \overline{S}$, by Theorem 1, $T$ is a critical set. This contradicts the assumption. ∎

Although a minimum perfect critical 3 set does not exist, a minimum perfect critical 4 set does exist. The topology structure of $G[S]$ will affect whether $S$ is an MPCS when $|S| \geq 4$.

### 3. MPCSs of $k$-distant Trees

*3.1 Main theorems*

Let $T$ be a tree. A *spine* in a tree $T$ is the longest path in $T$. A tree $T$ is $k$-distant if there is a spine $P$ such that every vertex of $T$ is the distance at most $k$ from a vertex in $P$. The paths are precisely 0-distant trees, 1-distant trees are exactly the caterpillars and the 2-distant trees are lobsters. Let spine $P = v_1 v_2 \cdots v_n$ be the longest path in $T$.

For path $P = v_1 v_2 \cdots v_n$, if the end node $v_1$ or $v_n$ be selected as the leader, path $P$ will be controllable. Hence, in this paper, we will discuss the leader set of lobsters with the structure shown in Figure 4.

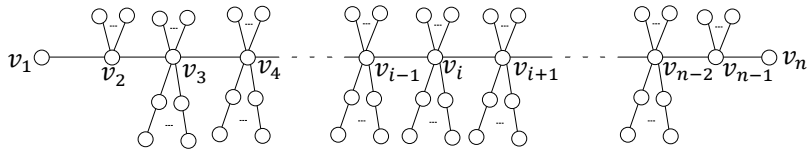

Figure 4  the structure of lobsters discussed in this paper

Let $p1(v_i), p2(v_i)$ be the number of vertices of distant 1, 2 from $v_i$, respectively. Let $S1(v_i)$, $S2(v_i)$ be the vertices set which contained all of the appendant vertices on the paths attached to $v_i$ with lengths 1, 2, respectively. By Theorem 2, any pair of vertices in $S1(v_i)$ is MPCS. Hence, in the rest of this paper, we set

$$|S \cap S1(v_i)| \leq 1 (\forall S \text{ is a } MPCS). \tag{3}$$

Similarly, we will prove in the following Theorem 4 that any four vertices which lie in a pair of paths adhered to $v_i$ form an MPCS. Hence, we set

$$|S \cap S2(v_i)| \leq 1. \tag{4}$$

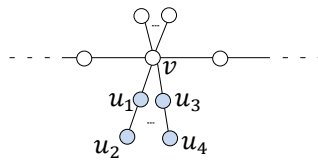

Figure 5　a lobster and a MPCS S=$\{u_1, u_2, u_3, u_4\}$

**Theorem 4**

Let $S = \{u_1, u_2, u_3, u_4\} \subset V(T)$ and $T$ is the *lobster* shown in Figure 5, then $S$ is a MPCS.

**Proof**  In this case, we have $\boldsymbol{L}_{S \to S} = \begin{bmatrix} 2 & -1 & 0 & 0 \\ -1 & 1 & 0 & 0 \\ 0 & 0 & 2 & -1 \\ 0 & 0 & -1 & 1 \end{bmatrix}$, $\lambda = \frac{3-\sqrt{5}}{2}$ is an eigenvalue of $\boldsymbol{L}_{S \to S}$, the corresponding eigenvector is $\boldsymbol{y}' = (x_1, x_2, x_3, x_4) = (1, \frac{\sqrt{5}+1}{2}, -1, -\frac{\sqrt{5}+1}{2})$. Notice that only one vertex, $v$, is adjacent to $S$, and $x_1 + x_3 = 0$, $S$ is a PCS of $T$. Further, by Theorem 2 and Theorem 3, any proper subset of $S$ is no longer a PCS. Hence, $S$ is a MPCS. ∎

The MPCSs obtained by Theorem 2 and Theorem 4 do not contain any vertices $v_i$ in spine path $P$. Next, the MPCS which include some $v_i \in V(P)$ will be discussed.

**Theorem 5**

Let $S$ be a MPCS, $P = v_1 v_2 \cdots v_n$ be the longest path in the lobster $T$. Let $v_{i_0} \in V(P) \cap S$ and $v_i \notin S(\forall\ 1 \le i \le i_0 - 1)$. If $p1(v_{i_0}) + p2(v_{i_0}) > 0$ and $d(v) \ge 3(\forall v \in S \cap V(P))$, then

　　(a) for $v_{i_0-1}$, $p2(v_{i_0-1}) \ge 1$.
　　(b) for $v_{i_0}$, $p1(v_{i_0}) = 1$ and $p2(v_{i_0}) = 0$.
　　(c) for $v_{i_0+1}$, $v_{i_0+1} \in S$.
　　(d) for $v_{i_0+2}$, $v_{i_0+2} \notin S$.

**Proof** (a) Since $v_{i_0} \in S$, $v_{i_0-1} \notin S$ and $v_{i_0-2} \notin S$, $p1(v_{i_0-1}) + p2(v_{i_0-1}) > 0$ can be derived from $N_S(v_{i_0-1}) \ne 1$(see Lemma 1).

Suppose $p2(v_{i_0-1}) = 0$, then $p1(v_{i_0-1}) \ge 1$. Let $u_1$ be the vertex adjacent to $v_{i_0-1}$ and $u_1 \in S$. Notice that $u_1$ is an appendant, $\lambda = 1$ is the corresponding eigenvalue for $S$. We claim that $S1(v_{i_0}) \cup S2(v_{i_0}) = \emptyset$. In fact, suppose $S1(v_{i_0}) \ne \emptyset$, let $u_2 \in S1(v_{i_0})$, then $u_2 \notin S$ since $\lambda = 1$ and $v_{i_0} \in S$. $N_S(u_2) = 1(u_2 \in \bar{S})$ is contradict to Lemma 1. Suppose $S2(v_{i_0}) \ne \emptyset$, from $|S \cap S2(v_{i_0})| \le 1(see\ (4))$, we have $|S \cap S2(v_{i_0})| = 1$. By Lemma 1, $d(v_{i_0}) = 3$. Let $u_3 \in S \cap S2(v_{i_0})$ and $u_4$ is the vertex connecting $v_{i_0}$ and $u_3$, then $u_4 \notin S$(see $\lambda = 1$). Let $S_0 = S \cap V(P)$, then we have

$$\boldsymbol{L}_{S_0 \to S_0} = \begin{bmatrix} 3 & -1 & 0 & 0 & \cdots & 0 & 0 \\ -1 & 3 & -1 & 0 & \cdots & 0 & 0 \\ 0 & -1 & 3 & -1 & \cdots & 0 & 0 \\ \vdots & \vdots & \vdots & \vdots & & \vdots & \vdots \\ 0 & 0 & 0 & 0 & \cdots & -1 & 3 \end{bmatrix},$$

Since $\boldsymbol{L}_{S_0 \to S_0} - I$ is a principal diagonally dominant matrix, $\lambda = 1$ isn't an eigenvalue of $\boldsymbol{L}_{S_0 \to S_0}$, it is a contradiction.

So far, we have $S1(v_{i_0}) \cup S2(v_{i_0}) = \emptyset$, it contradicts to the known condition $p1(v_{i_0}) + p2(v_{i_0}) > 0$.

(b) By (a), there exist a vertex, say $w_1$, such that $w_1$ is distant 2 from $v_{i_0-1}$ in path $P$. Let $w_2$ be the vertex which adjacent to $w_1$ and $v_{i_0-1}$. Then both $w_1$ and $w_2$ belongs to $S$. The corresponding eigenvalue $\lambda$ satisfy
$$(\lambda - 1)(\lambda - 2) = 1.$$

(5)

Suppose $p2(v_{i_0}) \neq 0$. Let $u_1 u_2 \in E(G)$ and $u_2 v_{i_0} \in E(G)$. If $u_2 \in \bar{S}$, then $u_1 \in S$, and $\lambda = 1$, this is a contradiction to (5). If $u_2 \in S$, then $u_1 \in S$. Consider the equation
$$L_{\{u_1,u_2\} \to V} y = \lambda y,$$
we have $\begin{cases} y_{u_1} - y_{u_2} = \lambda y_{u_1}, \\ -y_{u_1} + 2y_{u_2} - y_{v_{i_0}} = \lambda y_{u_2}. \end{cases}$ Notice that $(\lambda - 1)(\lambda - 2) = 1$, the equations means $y_{v_{i_0}} = 0$, it is a contradiction to $v_{i_0} \in S$ and $S$ is a MPCS.

Since $p2(v_{i_0}) = 0$ and $p1(v_{i_0}) + p2(v_{i_0}) > 0$, together with (3), we have $p1(v_{i_0}) = 1$.

(c) Suppose $v_{i_0+1} \notin S$, from (b), let $w_3 v_{i_0} \in E(G)$, Consider the equation
$$L_{\{w_3, v_{i_0}\} \to V} y = \lambda y,$$
we have $\begin{cases} y_{w_3} - y_{v_{i_0}} = \lambda y_{w_3}, \\ -y_{w_3} + 3y_{v_{i_0}} = \lambda y_{v_{i_0}}. \end{cases}$ This equation means $(1 - \lambda)(3 - \lambda) = 1$ and it contradicts to (5).

(d) Similarly, as what we have proved in (b), we know that there exists a vertex $w_4$ such that $w_4 v_{i_0+1} \in E(G)$ and $w_4 \in S$. For convenience, let $y_1 = y_{w_3}, y_2 = y_{v_{i_0}}, y_3 = y_{v_{i_0+1}}, y_4 = y_{w_4}, y_5 = y_{v_{i_0+2}}$. From the equation
$$L_{\{w_3, v_{i_0}, w_4, v_{i_0+1}\} \to V} y = \lambda y,$$
we have $\begin{cases} y_1 - y_2 = \lambda y_{1_3}, \text{①} \\ -y_1 + 3y_2 - y_3 = \lambda y_2, \text{②} \\ y_4 - y_3 = \lambda y_4, \text{③} \\ -y_4 - y_2 - y_5 + 3y_3 = \lambda y_3. \text{④} \end{cases}$ By ①, ② and (5), $y_2 = y_3$. From ③, we have $y_3 = (2 - \lambda)y_3$. From ④, we have $y_3 = (1 - \lambda)y_4$. From ⑤, we have $y_5 = (3 - \lambda)y_3 - y_2 - y_4 = 0$. This shows that $v_{i_0+2} \notin S$. ∎

Other than twins and 4 vertices MPCS in Figure 5, Theorem 5 provides some other types of MPCS for lobsters, see Figure 6. The black vertices in Figure 6 are MPCSs.

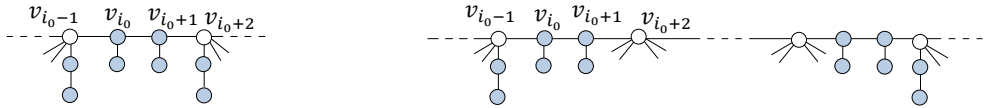

(a) MPCS with 8 vertices(filled)    (b)MPCS with $4n(n \geq 3)$ vertices(filled)

Figure 6   New MPCSs in *lobster*

*3.2 Algorithm*

From the above discussion, we have found four different types of MPCS. They are twins (Theorem 2), 4 vertices MPCS (Figure 5), 8 vertices MPCS (Figure 6(a)) and $4n(n \geq 3)$ vertices MPCS( Figure 6(b)). Based on these theories, this section will give an algorithm for finding the minimal leader set.

Critical Set Algorithm, CSA

Step0 input a lobster tree, set LeadersSet:=∅.

Step1 finds out the twins $S$ and take any $v \in S$, set LeadersSet:=LeadersSet ∪ {$v$}.

Step2 finds out 4 vertices MPCSs(Figure 5) $S$  and take any $v \in S$, set LeadersSet:=LeadersSet ∪ {$v$}.

Step3 If the lobster is controllable under LeadersSet, then output LeadersSet and Stop; otherwise, go to Step4.

Step4 finds out  $4n(n \geq 2)$ vertices MPCSs (Figure 6) $S$  and take any $v \in S$, set LeadersSet:=LeadersSet ∪ {$v$}.

Step5 If the lobster is controllable under LeadersSet, then output LeadersSet and Stop; otherwise, go to Step6.

Step6 find vertices $v_i$ such that $p1(v_{i-1}) + p2(v_{i-1}) > 0$ and $p1(v_i) + p2(v_i) = 0$, set LeadersSet:=LeadersSet $\cup \{v_i\}$. If the lobster is controllable under LeadersSet, then output LeadersSet and Stop; otherwise, output 'can't find'.

Altogether, there are four types of MPCSs provided above. Of course, there exist some other MPCSs. Hence, the last step of the algorithm is important. This step makes the algorithm more likely to find the leader set.

The algorithm flow chart for CSA is shown in Figure 7.

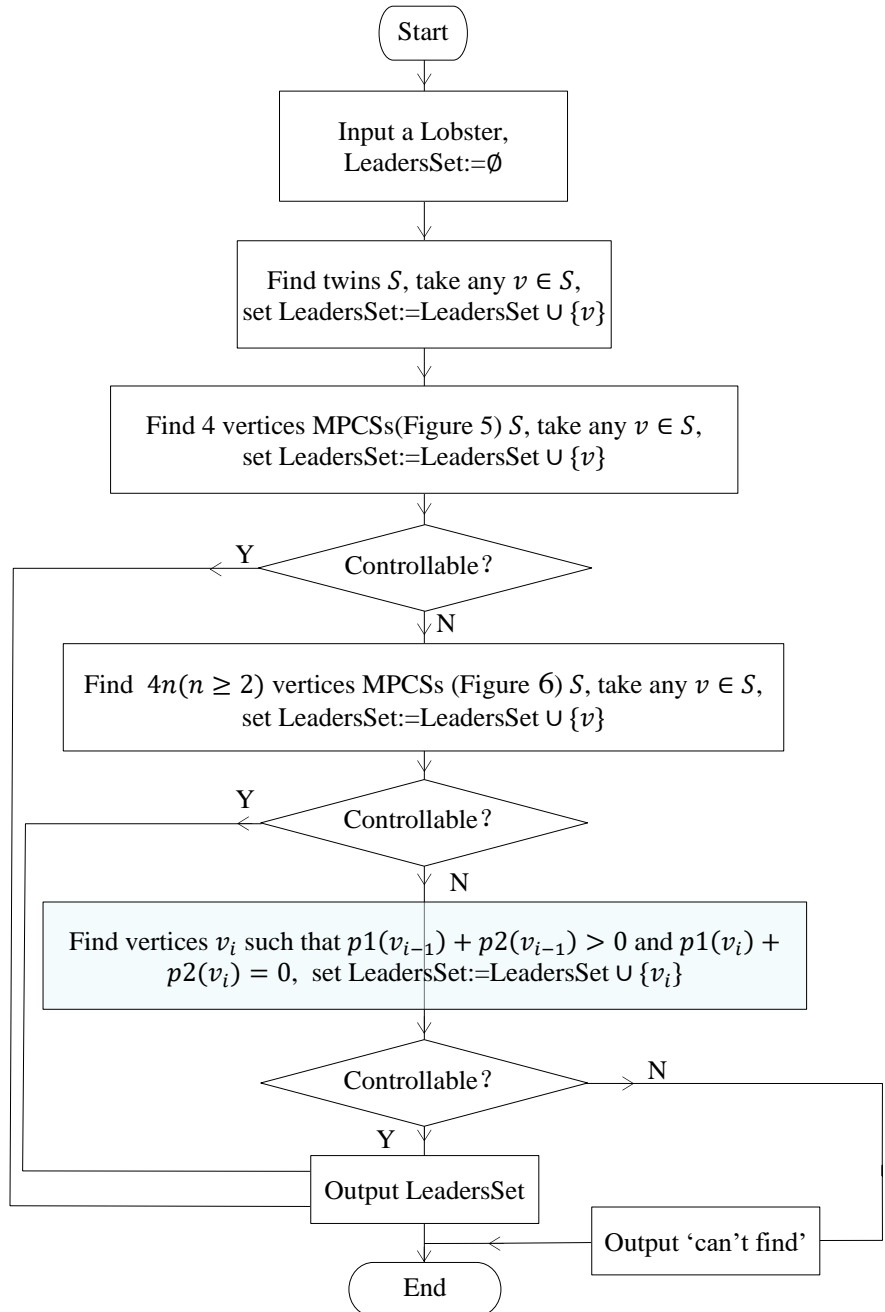

Figure 7  algorithm flow chart for finding leader set of lobster

## 4. Simulations
*4.1 Illustrate of the algorithm*

In this section, an example is given to illustrate the application of the theory and algorithm proposed. Graph $G$ in Figure 8 is a lobster. By Theorem 2, $S_1 = \{v_{10}, v_{11}\}$ and $S_2 = \{v_{16}, v_{17}\}$ are twins and MPCSs. By Theorem 4, $S_3 = \{v_1, v_2, v_{13}, v_{14}\}$ and $S_4 = \{v_{12}, v_{13}, v_{14}, v_{15}\}$ are MPCSs with 4 vertices. By Theorem 5,
$$S_5 = \{v_1, v_2, v_4, v_{18}, v_5, v_{19}, v_{20}, v_{21}\},$$
$$S_6 = \{v_{12}, v_{13}, v_4, v_{18}, v_5, v_{19}, v_{20}, v_{21}\},$$
$$S_7 = \{v_{14}, v_{15}, v_4, v_{18}, v_5, v_{19}, v_{20}, v_{21}\},$$
are MPCSs with 8 vertices. From Remark 1, leader set $F$ should satisfy that $F \cap S_i \neq \emptyset$. There are 320 such $F$'s. One of them is $F = \{v_4, v_{11}, v_{13}, v_{15}, v_{17}\}$(see filled vertices in Figure 8). Hence the minimum number of leaders in $G$ is 5. According to the theory about MPCS provided in this paper, the leader set can be found quickly. Otherwise, even if we know the minimum number of leaders, the probability of finding an $F$ is only 0.0157 by enumeration method.

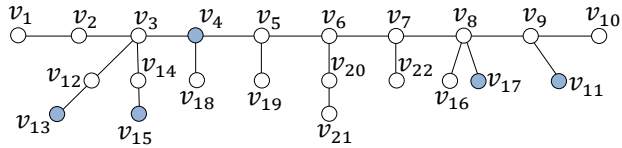

Figure 8  a lobster and its leader set

*4.2 Efficiency of the algorithm*

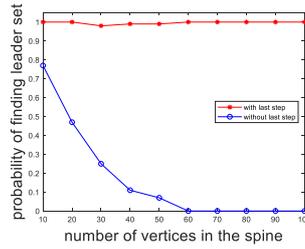

Figure 9  the probability of finding leader sets by the algorithm

In the experiment, let the spine $p = v_1 v_2 \cdots v_n$ and $10 \leq n \leq 100$. Note that $n$ is the number of vertices on the spine, not all vertices in the graph. Set $p1(v_i) + p2(v_i) \leq 2$. Paths of length 1 or 2 are randomly pasted on each vertex $v_i$. Figure 8 is such a randomly generated lobster. From the experimental results, the probability of finding leader sets is above 0.98. Through experiments, we also find that the last step (namely the step 'Find vertices $v_i$' in Figure 7) of the algorithm is important. Without this step, when the number of vertices is more than 50, the leader sets can hardly be found, see Figure 9.  This indicates that few MPCSs were found in the first four steps of the algorithm, and a large number of other MPCSs have not been found, which needs further study.

*4.3 Characteristics of leader set for lobster*

The following experiment is to explore the characteristics of a leader set of lobsters. It can be seen from Figure 10  Number of leaders  that the number of leaders increases with the number of vertices on the spine. And it increases linearly with $n$ approximately, that is  $l(n) \approx 0.3n + 2$.

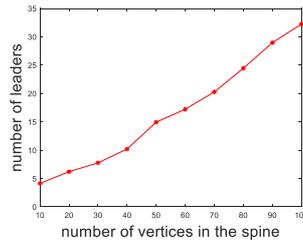

Figure 10   Number of leaders

Figure 11 studies the proportion of leader vertices. Let $N$ be the number of all vertices in the lobsters, $n$ is the number of vertices on the spine, and $l(n)$ is the number of leaders. See Figure 11, the proportion $\frac{l(n)}{N}$ is no more than 0.2. This shows that the random lobster graph is relatively easier to control.

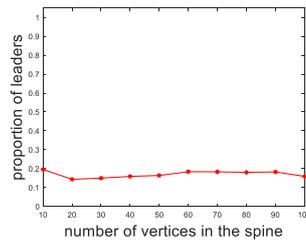

Figure 11   The proportion of leaders in the total number of vertices

## 5. Conclusion

In this paper, we study the controllability of the $k$-distant trees. The definition of MPCS is provided. Based on this new conception, an algorithm for finding the leaders in the lobsters is given, too. From the simulations, we see that the algorithm performs quite well, it can find out the leader set with a probability of more than 0.98. Further, some important characteristics of a leader set are discovered, such as the number of leaders and their proportion.

Although we have found four different types of MPCS in this paper, it is worth noting that a large number of MPCSs have not been found yet. Especially when $n \geq 50$, such MPCSs is likely to occur for the lobster graph.